\documentclass[11pt,a4paper]{article}
\usepackage{jheppub}
\usepackage{amsmath}
\usepackage[most]{tcolorbox}
\usepackage{dsfont}
\usepackage{ulem}
\usepackage{natbib}
\usepackage{xcolor}
\usepackage[hang,flushmargin]{footmisc}
\usepackage{tikz-cd}
\usepackage{enumitem}
\usepackage{amsfonts,amssymb, amscd,amsmath,latexsym,amsbsy,bm}
\usepackage{stmaryrd}
\usepackage{todonotes}

\usepackage{float}

\makeatletter\renewcommand{\@biblabel}[1]{#1.}\makeatother

\newtcolorbox{empheqboxed}{colback=gray!20, 
 colframe=white,
 width=\textwidth,
 sharpish corners,
 top=0mm, 
 bottom=0pt
}

\title{A remark on the q-hypergeometric integral Bailey pair and the solution to the star-triangle equation}

\author{Erdal Catak}
\affiliation{Department of Physics, Istanbul University,\\ 34134 Istanbul, Turkey \\[-0.4cm]

}

\emailAdd{ecatak@istanbul.edu.tr}

\abstract{We rewrite the recently constructed q-hypergeometric integral Bailey pair in a general form. Then with the help of the Bailey pair and $q$-beta hypergeometric sum-integral, we construct the star-triangle relation.}


\begin{document}
\maketitle
\flushbottom

\section{Introduction}
 Integral identities for trigonometric hypergeometric functions \cite{andrews1999special,schlosser2020hypergeometric, Spiridonov-essays, Catak:2021coz, Gahramanov:2015tta,Gahramanov:2016wxi,Bozkurt:2018xno,Mullahasanoglu:2021xyf,Rosengren:2016mnw,Hwang:2017kmk} have many applications in mathematical physics. Bailey lemma serves as a good tool for obtaining such complicated identities \cite{warnaar200150,spiridonov2002elliptic,spiridonov2004bailey}. Recently, several Bailey pairs were constructed via supersymmetric gauge theory computations, and integrable lattice spin models of statistical mechanics \cite{Gahramanov:2015cva,Brunner:2016uvv,Brunner:2017lhb,Spiridonov:2019uuw,Kashaev:2012cz,Gahramanov:2022jxz,Gahramanov:2021pgu}. 
 
 
 In \cite{Gahramanov:2015cva} the star-triangle relation \cite{baxter2016exactly,baxter:1997ssr} was constructed with help of integral Bailey pair by using the special condition that two consecutive pairs of variables are zero for the balancing condition of the $q$-beta sum-integral identity.  We reintroduce the sum-integral operator $M(t,m_t)_{x,m;z,n}$ and the $D$-function $D(t,m_t;x,m;z,n)$ in terms of $q$-hypergeometric function and rewrite Bailey pair. Namely, we show that it is possible to construct the Bailey pair  for  the $q$-beta sum-integral identity with a weaker balancing condition.

\section{Bailey pairs construction}

Here we consider the $q$-hypergeometric integral Bailey pair based on the $q$-hypergeometric integral-sum identity \cite{Gahramanov:2015cva, Kels:2015bda}.

A pair of sequences of functions $\beta_m(z;t,m_t)$ and $\alpha_m(z;t,m_t)$ is called a Bailey pair relative  to the parameter $t$, if the following integral-sum relation satisfied
\begin{align}
 \beta_m(x;t,m_t)=M(t,m_t)_{x,m;z,n}\alpha_n(z;t,m_t)\label{3.1}
\end{align}
where $|tx|,|t/x|<1$; $m,m_t,n\in\mathbb{Z}$; $x,z$ and $t\in\mathbb{C}$ and $M(t,m_t)_{x,m;z,n}$ is defined as
\begin{align}
    M(t,m_t)&_{x,m;x,n}:=\frac{(q^{m_t}t^2;q)_\infty}{(q^{1+m_t}t^{-2};q)_\infty}\sum_{n\in\mathbb{Z}}\int_{\mathbb{T}}[d_{n}z]
    \frac{1}{z^{2n}x^{2m}t^{2m_t}}\nonumber\\&\times\frac{(q^{1+\frac{m+m_{t}+n}{2}}\frac{1}{txz},q^{1+\frac{m+m_{t}-n}{2}}\frac{z}{tx};q)_\infty}{(q^{\frac{m+m_{t}+n}{2}}txz,q^{\frac{m+m_{t}-n}{2}}\frac{tx}{z};q)_\infty}\frac{(q^{1+\frac{-m+m_{t}+n}{2}}\frac{x}{tz},q^{1+\frac{-m+m_{t}-n}{2}}\frac{xz}{t};q)_\infty}{(q^{\frac{-m+m_{t}+n}{2}}\frac{tz}{x},q^{\frac{-m+m_{t}-n}{2}}\frac{t}{xz};q)_\infty}
\end{align}
with the following measure
\begin{align}
    [d_{m}z]:=\frac{(1-q^{m}z^2)(1-q^{m}z^{-2})}{q^m}\frac{dz}{4{\pi}iz}\;,\nonumber
\end{align}

and here $(z;q)_\infty=\prod_{j=0}^\infty(1-zq^j)$ is the $q$-Pochhammer symbol and the shorthand notation $(z;p,q)_\infty=(z;p)_\infty(z;q)_\infty$ is used. In order to construct a new Bailey pair we define the $D$-function as follows
\begin{align}
    D(t,n_t;\omega,n;z,m):&=\frac{q^{n_t}}{z^{2m}\omega^{2n}t^{2n_t}}\frac{(q^{1+\frac{n-n_{t}+m}{2}}q^{-\frac12}\frac{t}{z\omega},q^{1+\frac{n-n_{t}-m}{2}}q^{-\frac12}\frac{tz}{\omega};q)_\infty}{(q^{\frac{n-n_{t}+m}{2}}q^{\frac12}\frac{\omega{z}}{t},q^{\frac{n-n_{t}-m}{2}}q^{\frac12}\frac{\omega}{tz};q)_\infty}\nonumber\\&\times\frac{(q^{1+\frac{-n-n_{t}+m}{2}}q^{-\frac12}\frac{t\omega}{z},q^{1+\frac{-n-n_{t}-m}{2}}q^{-\frac12}t\omega{z};q)_\infty}{(q^{\frac{-n-n_{t}+m}{2}}q^{\frac12}\frac{z}{t\omega},q^{\frac{n-n_{t}-m}{2}}q^{\frac12}\frac{1}{t\omega{z}};q)_\infty} .
\end{align}
This function satisfies the following properties: $D(t^{-1},-n_t;\omega,n;z,m)D(t,n_t;\omega,n;z,m)=1$ and $D(1,0;\omega,n;z,m)=1$. These properties are important for the construction of the so-called star-triangle relation. 

Let $\beta_m(x;t,m_t)$ and $\alpha_m(x;t,m_t)$ be a Bailey pair then  
\begin{align} 
\alpha'_{k}(x;st,n_s,n_t)&=D(s,n_s;y,l;x,k)\alpha_{k}(x;t,n_t)\label{3.4},\\ \nonumber
\beta'_{k}(x;st,n_s,n_t)&=D(t^{-1},n_t;y,l;x,k)M(s,n_s)_{x,k;z,m}D(st,n_s,n_t;y,l;z,m)\\ 
& \qquad \times\beta_{m}(z;t,n_t), \label{3.5}
\end{align}
also form a Bailey pair via  the Bailey lemma \cite{Magadov:2018hlc,spiridonov2004bailey,Gahramanov:2015cva}. Here $s,y\in\mathbb{C}$ and $l,n_{t,s}\in\mathbb{Z}$  \cite{Magadov:2018hlc}. From (\ref{3.1}), (\ref{3.4}) and (\ref{3.5}) one finds the following star-triangle relation
\begin{align}
    M(s,n_s)&_{\omega,k;z,m}D(st,n_s,n_t;y,l;z,m)M(t,n_t)_{z,m;x,j}\alpha_j(x;t,n_t)\nonumber\\&=D(t,n_t;y,l;\omega,k)M(st,n_s,n_t)_{\omega,k;x,j}D(s,n_s;y,l;x,j)\alpha_j(x;t,n_t) \;\label{3.6}
\end{align}
This relation is crucial for the integrability of the corresponding integrable system. The equality (\ref{3.6}) was obtained in the \cite{Gahramanov:2015cva} with the help of the $q$-beta sum-integral \cite{Gahramanov:2016ilb,Rosengren:2016mnw} 
\begin{align}
    \sum_{m\in\mathbb{Z}}\int_{\mathbb{T}}\prod_{j=1}^{6}[d_{m}z]\frac{(q^{1+\frac{m+n_j}{2}}\frac{1}{a_{j}z},q^{1+\frac{n_j-m}{2}}\frac{z}{a_{j}};q)_\infty}{z^{6m}(q^\frac{m+n_j}{2}a_{j}z,q^\frac{n_j-m}{2}\frac{a_{j}}{z};q)_\infty}=\frac{1}{\prod_{j=1}^{6}a_j^{n_j}}\prod_{1\leq j<k\leq 6}\frac{(q^{1+\frac{n_j+n_k}{2}}\frac{1}{a_{j}a_{k}};q)_\infty}{(q^{\frac{n_j+n_k}{2}}a_{j}a_{k};q)_\infty} \;,\label{2.7}
\end{align}
with the balancing conditions
\begin{equation} \label{balanc}
\sum_{j=1}^6n_j=0 \;\; \text{and} \;\; \prod_{j=1}^6a_j=q.
\end{equation}
where $a_j,q\in\mathbb{C},$ and $n_j\in\mathbb{Z}$, $j=1,...,6$. 

In  \cite{Gahramanov:2015cva} for the construction of the Bailey pair authors  use stronger balancing condition for the integral identity (\ref{2.7}), namely the sum of successive pairs is zero, i.e. $n_1+n_2=0$, $n_3+n_4=0$ and $n_5+n_6=0$. In this work we show that the (\ref{3.6}) holds for any $a_j,q\in\mathbb{C},$ and $n_j\in\mathbb{Z}$, $j=1,...,6$ satisfying the balancing conditions (\ref{balanc}).

It is sufficient to compute the left-hand side of the expression (\ref{3.6}) with  the help of $q$-beta hypergeometric integral-sum identity (\ref{2.7}) 
\begin{align}
     M&(s,n_s)_{\omega,k;z,m}D(st,n_s,n_t;y,l;z,m)M(t,n_t)_{z,m;x,j}\nonumber
  \\
  =&\frac{(q^{n_s}s^2,q^{n_t}t^2;q)_\infty}{(q^{1+n_s}s^{-2},q^{1+n_t}t^{-2};q)_\infty}\sum_{m\in\mathbb{Z}}\int_{\mathbb{T}} [d_{m}z]\frac{1}{\omega^{2k}z^{2m}s^{2n_s}}\frac{(q^{1+\frac{k+n_{s}+m}{2}}\frac{1}{s\omega{z}},q^{1+\frac{k+n_{s}-m}{2}}\frac{z}{s\omega};q)_\infty}{(q^{\frac{k+n_{s}+m}{2}}s\omega{z},q^{\frac{k+n_{s}-m}{2}}\frac{s\omega}{z};q)_\infty}\nonumber
  \\
  &\times\frac{(q^{1+\frac{-k+n_{s}+m}{2}}\frac{\omega}{sz},q^{1+\frac{-s+n_{s}-m}{2}}\frac{\omega{z}}{s};q)_\infty}{(q^{\frac{-k+n_{s}+m}{2}}\frac{sz}{\omega},q^{\frac{-k+n_{s}-m}{2}}\frac{s}{\omega{z}};q)_\infty}\nonumber
  \\
  &\times\frac{q^{n_t+n_s}}{z^{2m}y^{2l}t^{2n_t+2n_s}s^{2n_s+2_t}}\frac{(q^{1+\frac{k-n_{s}-n_{t}+m}{2}}q^{-\frac12}\frac{ts}{yz},q^{1+\frac{k-n_{s}-n_{t}-m}{2}}q^{-\frac12}\frac{tsz}{y};q)_\infty}{(q^{\frac{k-n_{s}-n_{t}+m}{2}}q^{\frac12}\frac{yz}{ts},q^{\frac{k-n_{s}-n_{t}-m}{2}}q^{\frac12}\frac{y}{tsz};q)_\infty}\nonumber
  \\
  &\times\frac{(q^{1+\frac{-k-n_{s}-n_{t}+m}{2}}q^{-\frac12}\frac{tsy}{z},q^{1+\frac{-k-n_{s}-n_{t}-m}{2}}q^{-\frac12}tsyz;q)_\infty}{(q^{\frac{-k-n_{s}-n_{t}+m}{2}}q^{\frac12}\frac{z}{tsy},q^{\frac{-k-n_{s}-n_{t}-m}{2}}q^{\frac12}\frac{1}{styz};q)_\infty}\nonumber 
  \\
  &\times\sum_{j\in\mathbb{Z}}\int_{\mathbb{T}}[d_{j}x]\frac{1}{z^{2m}x^{2j}t^{2n_t}}\frac{(q^{1+\frac{m+n_{t}+j}{2}}\frac{1}{txz},q^{1+\frac{m+n_{t}-j}{2}}\frac{x}{tz};q)_\infty}{(q^{\frac{m+n_{t}+j}{2}}txz,q^{\frac{m+n_{t}-j}{2}}\frac{tz}{x};q)_\infty}\frac{(q^{1+\frac{-m+n_{t}+j}{2}}\frac{z}{tx},q^{1+\frac{-m+n_{t}-j}{2}}\frac{xz}{t};q)_\infty}{(q^{\frac{-m+n_{t}+j}{2}}\frac{tx}{x},q^{\frac{-m+n_{t}-j}{2}}\frac{t}{xz};q)_\infty}\nonumber
  \\
  =&\frac{(q^{n_s}s^2,q^{n_t}t^2;q)_\infty}{(q^{1+n_s}s^{-2},q^{1+n_t}t^{-2};q)_\infty}\sum_{j\in\mathbb{Z}}\int_{\mathbb{T}}[d_{j}x]\frac{q^{n_t+n_s}}{\omega^{2k}x^{2j}y^{2l}s^{4n_s+2n_t}t^{4n_t+2n_s}}\nonumber
  \\
  &\times\sum_{m\in\mathbb{Z}}\int_{\mathbb{T}}[d_{m}z]\prod_{j=1}^{6}\frac{(q^{1+\frac{m+n_j}{2}}\frac{1}{a_{j}z},q^{1+\frac{n_j-m}{2}}\frac{z}{a_{j}};q)_\infty}{z^{6m}(q^\frac{m+n_j}{2}a_{j}z,q^\frac{n_j-m}{2}\frac{a_{j}}{z};q)_\infty}\label{2.9},
\end{align}
where we introduced new parameters as follows
\begin{align}
    \begin{aligned}
         a_1&=s\omega,  &n_1&=k+n_s, &a_4&=q^{\frac12}(st)^{-1}y^{-1},  &n_4&=-l-n_s-n_t,
         \\
         a_2&=s\omega^{-1},  &n_2&=-k+n_s,,    &a_5&=tx,  &n_5&=j+n_t,
         \\
         a_3&=q^{\frac12}(st)^{-1}y,  &n_3&=l-n_s-n_t,  &a_6&=tx^{-1},  &n_6&=-j+n_t\label{2.10},
    \end{aligned}
\end{align}
and one can check that (\ref{2.10}) the balancing conditions are satisfied.

Now we use the integral identity (\ref{2.7}) and end up  with the following expression
\begin{align}
\frac{(q^{n_s}s^2,q^{n_t}t^2;q)_\infty}{(q^{1+n_s}s^{-2},q^{1+n_t}t^{-2};q)_\infty}\sum_{j\in\mathbb{Z}}\int_{\mathbb{T}}[d_{j}x]\frac{q^{n_t+n_s}}{\omega^{2k}x^{2j}y^{2l}s^{(4n_s+2n_t)}t^{(4n_t+2n_s)}}\frac{1}{\prod_{j=1}^{6}a_j^{n_j}}\prod_{1\leq j<r\leq 6}\frac{(q^{1+\frac{n_j+n_r}{2}}\frac{1}{a_{j}a_{r}};q)_\infty}{(q^{\frac{n_j+n_r}{2}}a_{j}a_{r};q)_\infty}\label{2ff}.
\end{align}
Rearrangement of expression (\ref{2ff}) and using following identity \cite{Dimofte:2011py,Gahramanov:2016wxi,Gahramanov:2015cva}
\begin{align}
    \frac{(q^{1-n_s-n_t}(st)^{-2};q)_\infty}{(q^{-n_s-n_t}(st)^2;q)_\infty}=\frac{q^{n_s+n_t}}{(-(st)^2)^{2(n_s+n_t)}}\frac{(q^{1+n_s+n_t}(st)^{-2};q)_\infty}{(q^{n_s+n_t}(st)^2;q)_\infty},    \   \  \    n_s,n_t
    \in\mathbb{Z},
\end{align}
yields the right-hand side of the expression (\ref{3.6})
\begin{align}
    \begin{aligned}
         \frac{q^{n_t}(q^{1+\frac{k+l-n_t}{2}}q^{-\frac12}\frac{t}{y\omega},q^{1+\frac{-k+l-n_t}{2}}q^{-\frac12}\frac{t\omega}{y},q^{1+\frac{k-l-n_t}{2}}q^{-\frac12}\frac{ty}{\omega},q^{1+\frac{-k-l-n_t}{2}}q^{-\frac12}ty\omega;q)_\infty}{\omega^{2k}y^{2l}t^{2n_t}(q^{\frac{k+l-n_t}{2}}q^{\frac12}\frac{y\omega}{t},q^{\frac{-k+l-n_t}{2}}q^{\frac12}\frac{y}{t\omega},q^{\frac{k-l-n_t}{2}}q^{\frac12}\frac{\omega}{ty},q^{\frac{-k-l-n_t}{2}}q^{\frac12}\frac{1}{ty\omega};q)_\infty}
         \\
         \times\frac{q^{n_s+n_t}}{(st)^{4(n_s+n_t)}}\left(\frac{(st)^{4(n_t+n_s)}}{q^{n_s+n_t}}\frac{(q^{n_s+n_t}(st)^2;q)_\infty}{(q^{1+n_s+n_t}(st)^{-2};q)_\infty}\right)\sum_{j\in\mathbb{Z}}\int_{\mathbb{T}}[d_{j}x]\frac{1}{\omega^{2k}x^{2j}t^{2n_t}s^{2n_s}}
         \\
         \times\frac{(q^{1+\frac{k+j+n_s+n_t}{2}}\frac{1}{stx\omega},q^{1+\frac{k-j+n_s+n_t}{2}}\frac{x}{st\omega},q^{1+\frac{-k+j+n_s+n_t}{2}}\frac{\omega}{stx},q^{1+\frac{-k-j+n_s+n_t}{2}}\frac{x\omega}{st};q)_\infty}{(q^{\frac{k+j+n_s+n_t}{2}}stx\omega,q^{\frac{k-j+n_s+n_t}{2}}\frac{st\omega}{x},q^{\frac{-k+j+n_s+n_t}{2}}\frac{stx}{\omega},q^{\frac{-k-j+n_s+n_t}{2}}\frac{st}{x\omega};q)_\infty}
         \\
         \times\frac{q^{n_s}(q^{1+\frac{l+j-n_s}{2}}q^{-\frac12}\frac{s}{yx},q^{1+\frac{l-j-n_s}{2}}q^{-\frac12}\frac{sx}{y},q^{1+\frac{-l+j-n_s}{2}}q^{-\frac12}\frac{sy}{x},q^{1+\frac{-l-j-n_s}{2}}q^{-\frac12}syx;q)_\infty}{y^{2l}x^{2j}s^{2n_s}(q^{\frac{l+j-n_s}{2}}q^{\frac12}\frac{yx}{s},q^{\frac{l-j-n_s}{2}}q^{\frac12}\frac{y}{sx},q^{\frac{-l+j-n_s}{2}}q^{\frac12}\frac{x}{sy},q^{\frac{-l-j-n_s}{2}}q^{\frac12}\frac{1}{syx};q)_\infty}
         \\
=D(t,n_t;y,l;\omega,k)M(st,n_s,n_t)_{\omega,k;x,j}D(s,n_s;y,l;x,j).
    \end{aligned}
\end{align}
As one can see, we generalized $D$ function and operator $M$ (and the Bailey pair respectively) and obtained the star-triangle relation discussed in \cite{Gahramanov:2015cva, Kels:2015bda}.

\section{Conclusion}

To conclude, we show that it is not necessary to use special conditions for balancing conditions for $q$-beta hypergeometric sum-integral identity to prove the corresponding star-triangle relation via the Bailey pair construction. Using our result one can apply Bailey pairs sequentially to obtain more complicated $q$-hypergeometric sum-integral identities, see e.g. \cite{spiridonov2004bailey,Mullahasanoglu:2021xyf}.  


\section*{Acknowledgements}
We would like to thank Ilmar Gahramanov for providing us with the problem and for his many useful comments. We are also grateful to Mustafa Mullahasanoğlu for their valuable discussions and to Boğaziçi University for hospitality where part of this work was done.


\appendix


\bibliographystyle{utphys}
\bibliography{refYBE}



\end{document}